\documentclass[11pt,reqno]{amsart}
\usepackage{amsmath,amssymb,amsthm,enumerate,natbib,color,bbm,ifthen,graphicx,overpic}
\usepackage[latin1]{inputenc}

\def\nset{{\mathbb{N}}}
\def\rset{\mathbb R}

\newcommand{\eps}{\varepsilon}

\def\F{\mathcal{F}} 
\def\B{\mathcal{B}} 

\def\E{\mathbb{E}}

\def\M{\mathcal{M}}

\def\T{\mathbb{T}}
\def\A{\mathcal{A}}

 \newcommand{\X}{\mathcal{X}}

\def\PE{\mathbb{E}} 

\def\L{\mathcal{L}} 

\def\barV{\overline{V}}
\def\barW{\overline{W}}
\def\v{\textsf{v}}

\def\C{\textsf{C}} 

\newcommand{\nnorm}[1]{\left\vert\!\left\vert\!\left\vert#1\right\vert\!\right\vert\!\right\vert}

\usepackage{bm}

\newlength{\noteWidth}
\newtheorem{theo}{Theorem}[section]
\newtheorem{lemma}[theo]{Lemma}

\newtheorem{prop}[theo]{Proposition}
\newtheorem{defi}[theo]{Definition}

\theoremstyle{remark}
\newtheorem{rem}{Remark}

\newcounter{hypoconbis}
\newcounter{saveconbis}
\newcommand\debutA{\begin{list} {\textbf{Assumption A\arabic{hypoconbis}}}{\usecounter{hypoconbis}}\setcounter{hypoconbis}{\value{saveconbis}}}
\newcommand\finA{\end{list}\setcounter{saveconbis}{\value{hypoconbis}}}

\newcounter{hypoconbisp}
\newcounter{saveconbisp}
\newcommand\debutAp{\begin{list} {\textbf{Assumption A\arabic{hypoconbisp}'}}{\usecounter{hypoconbisp}}\setcounter{hypoconbisp}{\value{saveconbisp}}}
\newcommand\finAp{\end{list}\setcounter{saveconbisp}{\value{hypoconbisp}}}

\newcounter{hypocom}
\newcounter{savecom}
\newcommand{\debutB}{\begin{list}{\textbf{B\arabic{hypocom}}}{\usecounter{hypocom}}\setcounter{hypocom}{\value{savecom}}}
\newcommand{\finB}{\end{list}\setcounter{savecom}{\value{hypocom}}}

\newcounter{hypocomp}
\newcounter{savecomp}
\newcommand{\debutBp}{\begin{list}{\textbf{B\arabic{hypocomp}'}}{\usecounter{hypocomp}}\setcounter{hypocomp}{\value{savecomp}}}
\newcommand{\finBp}{\end{list}\setcounter{savecomp}{\value{hypocomp}}}

\newcounter{hypostab}
\newcounter{savestab}
\newcommand{\debutC}{\begin{list}{\textbf{C\arabic{hypostab}}}{\usecounter{hypostab}}\setcounter{hypostab}{\value{savestab}}}
\newcommand{\finC}{\end{list}\setcounter{savestab}{\value{hypostab}}}

\newcounter{hypodist}
\newcounter{savedist}
\newcommand{\debutD}{\begin{list}{\textbf{D\arabic{hypodist}}}{\usecounter{hypodist}}\setcounter{hypodist}{\value{savedist}}}
\newcommand{\finD}{\end{list}\setcounter{savedist}{\value{hypodist}}}

\begin{document}

\title[Quadratic forms of Markov chains]{Limit theorems for quadratic forms of Markov Chains}

\author{Yves F. Atchad\'e} \thanks{ Y. F. Atchad\'e: University of Michigan, Department of Statistics, Ann Arbor,
 48109, MI, United States. {\em E-mail address:} yvesa@umich.edu}
\author{Matias D. Cattaneo} \thanks{ Matias D. Cattaneo: University of Michigan, Department of Economics, Ann Arbor,
 48109, MI, United States. {\em E-mail address:} cattaneo@umich.edu}



\maketitle

\begin{center} (April 2011) \end{center}

\begin{abstract}
We develop a martingale approximation approach to studying the limiting behavior of quadratic forms of Markov chains. 
We use the technique to examine the asymptotic behavior of lag-window estimators in time series and  we apply the results to Markov Chain Monte Carlo simulation.  As another illustration, we use the  method to derive a central limit theorem for U-statistics with varying kernels.
\end{abstract}

\bigskip

\setcounter{secnumdepth}{3}

\section{Introduction}
This paper deals with quadratic forms of the type 
\begin{equation}\label{Un}
U_n(h_n)=\sum_{\ell=1}^n\sum_{j=1}^\ell w_n(\ell,j)h_n(X_\ell,X_j),\qquad n\geq 1,\end{equation}
 for a stochastic process $\{X_n,\;n\geq 0\}$, weight matrices $w_n:\;\nset\times\nset\to\rset$ and symmetric kernels $h_n:\;\X\times\X\to\rset$. Quadratic forms of possibly time-dependent random variables naturally arise in a variety of statistical and econometric problems, and their asymptotic properties are of particular importance to develop asymptotically valid inference procedures.

For independent sequences $\{X_n,\;n\geq 0\}$, the well known Hoeffding decomposition provides a useful approach to studying the asymptotic properties of $U_n(h_n)$ because it decomposes the statistic into two (uncorrelated) martingale sequences, which are then easily handled by standard martingale theory. See, e.g., \cite{Serfling_1980_book} for a review. When the process $\{X_n,\;n\geq 0\}$ is time-dependent, however, the classical Hoeffding decomposition is not very useful because the resulting representation does not have the desirable martingale property in general. As a consequence, the large sample properties of quadratic forms of time-dependent random variables are typically established in a less systematic way. The most well understood case is the case of a standard U-statistics where $h_n$ does not depend on $n$ and $w_n(\ell,j)=1$ if $\ell\neq j$ and $0$ otherwise (\cite{yoshihara76,eagleson79,dehlingetal10}). There has been some recent progress. \cite{hsingetwu04} considers $U_n(h)$ where neither $h_n$ nor $w_n$ depends on $n$, whereas \cite{wuetshao07} studies $U_n(h_n)$ when $h_n(x,y)=h(x,y)=xy$ for a martingale-difference sequence (see also \cite{bhansalietal07} for i.i.d. sequences).

We develops a martingale approximation for $U_n(h_n)$ which allows for a general and systematic analysis of $U_n(h_n)$ when $\{X_n,\;n\geq 0\}$ is a Markov chain. Martingale approximation is a well established technique when dealing with linear partial sums of dependent processes (\cite{mw00,merlevedeetal06}), but has not been fully explored in dealing with quadratic forms (a notable exception is \cite{wuetshao07}). In the present paper, we obtain an approximating quadratic martingale to $U_n(h_n)$ from a solution of a bivariate analog of the well known Poisson's equation. 

As an application we  study the asymptotic behavior of lag-window estimators of long-run variance (asymptotic variance)  for Markov chains (see, e.g., \cite{priestley81}).  We obtain a decomposition of lag-window estimators that shed some new light on the asymptotic behavior of these estimators, particularly by contrasting the classical asymptotics and the so-called "fixed-b" asymptotics (\cite{neave70,kieferetvogelsang05}). We derive two theorems that extend existing results. We obtain the consistency of lag-window estimators for non-geometrically ergodic Markov chains extending recent results of \cite{jonesetal09} and \cite{atchade10}; and we extend the "fixed-b" asymptotics framework to handle non-stationary Markov chains. These results have important implications for Markov Chain Monte Carlo (MCMC) simulations, offering in particular new robust procedures for constructing Monte Carlo confidence intervals.

As another application of the martingale approximation method, we derive a central limit theorem for U-statistics with varying kernels without imposing stationarity and under assumptions that are more easily verifiable. In particular, we do not rely on mixing conditions.

The paper is organized as follows. The rest of the introduction outlines the general setup and introduces the main notation employed throughout, while Section 2 derives the main martingale approximation method. Section 3 derives the asymptotic properties of lag-window estimators and, in particular, applies these results to MCMC simulation.  We study U-statistics with varying kernels in Section \ref{sec:Ustat}. All the proofs are presented in Section \ref{sec:proofs}.

\subsection{Setup and Notation}
Throughout the paper, $\{X_n,\;n\geq 0\}$ denotes a Markov chain taking values in a general state space $(\X,\B)$ equipped with a countably generated sigma-algebra $\B$. We denote by $P$ the transition kernel of the Markov chain and $\mu$ its invariant distribution whose existence is assumed. Unless explicitly stated otherwise, $\{X_n,\;n\geq 0\}$ is a nonstationary Markov chain with initial distribution $\rho$.

We will rely on the following set of general notation. Suppose that $(\T,\A)$ be an arbitrary measure space. If $W: \mathbb{T}\to [1, +\infty)$ is a function, the $W$-norm of a function
$f: \mathbb{T}\to \rset$ is defined as $|f|_W:=\sup_{x\in\mathbb{T}} |f(x)|
/W(x)$. The set of measurable functions $f:\;\T\to\rset$ with finite
$W$-norm is denoted by $\L_W(\T)$ or simply $\L_W$ when there is no ambiguity on the space $\T$. For a finite real-valued signed measure $\nu$ on $\mathbb{T}$, we denote the $W$-norm of $\nu$ as
\[\|\nu\|_{W}:=\int W(x)|\nu|(dx)=\sup_{|f|_{W}\leq 1}\left|\int f(x)\nu(dx)\right|,\]
where $|\nu|$ is the total variation measure of $\nu$. We denote $\M_W(\mathbb{T})$ the space of all finite real-valued signed measures $\nu$ on $\mathbb{T}$ such that $\|\nu\|_W<\infty$. It is well-known that $(\M_W(\mathbb{T}),\|\cdot\|_W)$ is a Banach space. When the measure space $\mathbb{T}$ is understood, we simply write $\M_W$. We will use the notation $\nu(f)$ to denote the integral $\int f(x)\nu(dx)$.  If $\mu,\nu$ are two finite signed measures on $(\mathbb{T},\A)$, we denote their product by $\mu\nu$ or $\mu\bigotimes\nu$, and the product of a finite number $k$ of finite signed measures $\nu_1,\ldots,\nu_k$ is denoted by $\bigotimes_{j=1}^k\nu_j$. 

If $Q$ is a transition kernel on $(\T,\A)$, its iterates are defined as: $Q^0$ is the identity kernel ($Q^0(x,A)=\textbf{1}_A(x)$) and for $n\geq 1$, we define $Q^n(x,\cdot)=\int Q(x,dz)Q^{n-1}(z,\cdot)$. If $h:\;\T\times\T\to\rset$ is a bivariate function then $Qh$ is the bivariate function defined by the rule $Qh(x,y)=\int Q(x,dz)h(z,y)$ and $Q^2h$ is defined as $Q^2h(x_1,x_2)=\int Q(x_1,dz_1)\int Q(x_2,dz_2)h(z_1,z_2)$. If $h:\;\T\to\rset$ is univariate, $Qh$ is defined similarly as $Qh(x)=\int Q(x,dz)h(z)$. Fix $Q$ a Markov kernel, and $V:\;\T\times \T\to [1,\infty)$. For $p\geq 1$ and a function $h:\;\T\times\T\to\rset$, we define
\[ \nnorm{h}_{p,V}:=\sup_{x,y\in \T} \frac{\left(\int Q(x,dz)|h(z,y)|^p\right)^{1/p}}{V(x,y)}.\]
For a univariate function $V:\; \T\to [1,\infty)$ and for $h:\;\T\to\rset$, we define $\nnorm{h}_{p,V}$ similarly as 
\[\nnorm{h}_{p,V}:=\sup_{x\in \T} V(x)^{-1}\left(\int Q(x,dz)|h(z)|^p\right)^{1/p}.\] 
When we use the notation $\nnorm{h}_{p,V}$ below, it will always be with respect to $P$, the Markov kernel of the reference process $\{X_n,\;n\geq0\}$, unless stated otherwise. 
The following short-range dependence concept will play an important rule.
\begin{defi}\label{def1}
Fix $r\in\nset$. For measurable functions $\barV_r\leq \barW_r:\T^r\to[1,\infty)$, we say that the transition kernel $Q$ with invariant distribution $\mu$ satisfies the condition $\C(r,\barV_r,\barW_r)$ if there exists a finite constant $c$ such that 
\begin{equation}\label{ergocond}
\sum_{\ell_1\geq 0}\cdots\sum_{\ell_r\geq 0}\left\|\bigotimes_{j=1}^r\left(Q^{\ell_j}(x_j,\cdot)-\mu\right)\right\|_{\barV_r}\leq c \barW_r(x_1,\ldots,x_r),\;\;(x_1,\ldots,x_r)\in \X^r.\end{equation}
\end{defi}
%

Throughout the paper, we denote by $c$ a finite constant which depends solely on the kernel $P$ but whose actual value can change from one equation to the next. In particular $c$ does not depend on the family of function $\{h_n,\;n\geq 1\}$ considered. Finally, all limits are taken as $n \to \infty$ unless explicitly noted otherwise.

\section{A martingale approximation for quadratic forms}
For notational convenience, we shall write $\bar\mu$ to denote the product probability measure $\bar\mu(du,dv)=\mu(du)\mu(dv)$, where $\mu$ is the invariant distribution of the Markov kernel $P$.  Consider the following assumption.
\vspace{0.2cm}
\debutA
\item \label{A1} There exist symmetric measurable functions $\barV_2\leq \barW_2:\;\X\times\X\to [1,\infty)$ such that $P$ satisfies $\C(2,\barV_2,\barW_2)$. Furthermore, $P^s\barW_2(x)<\infty$ for all $x\in\X^2$ and for $s\in\{1,2\}$. 
\finA

\begin{rem}
It is always possible to deduce A\ref{A1} from a univariate short-range dependence assumption.  Indeed, if $P$ satisfies $\C(1,V_1,W_1)$ and $\C(1,V_2,W_2)$, and $PW_1<\infty$, $PW_2<\infty$, define $\barV_2(x,y)=V_1(x)V_2(y)$ and $\barW_2(x,y)=W_1(x)W_2(y)$. Then 
\[\|\left(P^n(x,\cdot)-\mu\right)\bigotimes\left(P^m(y,\cdot)-\mu\right)\|_{\barV_2}=\|P^n(x,\cdot)-\mu\|_{V_1}\|P^m(y,\cdot)-\mu\|_{V_2}.\]
Thus
\[\sum_{n\geq 0}\sum_{m\geq 0}\|\left(P^n(x,\cdot)-\mu\right)\bigotimes\left(P^m(y,\cdot)-\mu\right)\|_{\barV_2}\leq c\barW_2(x,y),\]
and therefore A\ref{A1} holds.
\end{rem} 

\begin{rem}\label{remA1}
The univariate  condition $\C(1,V,W)$ holds for geometrically ergodic Markov kernels (that is, kernels $P$ for which $\left\|P^n(x,\cdot)-\mu\right\|_{V}$ converges to zero exponentially fast for some $V\geq 1$). It also holds for sub-geometrically ergodic Markov kernels ($\left\|P^n(x,\cdot)-\mu\right\|_{V}$ converges to zero sub-geometrically) for which the rate of convergence is summable. It is sometimes possible to check the condition $\C(1,V,W)$ using Lyapunov drift conditions and their extensions and this has been done for several time series Markov models (\cite{doucetal04,meitzetsaikkonen08,meynettweedie93}). 
\end{rem}

We show that whenever A\ref{A1} holds, there exists a martingale approximation to $U_n(h_n)$ that offers a simple route to study the asymptotics of $U_n(h_n)$.  The space $\M_{\barV_2}(\X\times\X)$ of all finite signed measure on $\X\times\X$ with finite $\|\cdot\|_{\barV_2}$ norm, equipped with the norm $\|\cdot\|_{\barV_2}$ is a Banach space. Under A\ref{A1} and for any $x,y\in\X$, 
\[\bar R_2(x,y;(du,dv)):=\sum_{n_1\geq 0}\sum_{n_2\geq 0}\left(P^{n_1}(x,du)-\mu(du)\right)\bigotimes\left(P^{n_2}(y,dv)-\mu(dv)\right)\]
is a finite signed measure that  belongs to $\M_{\barV_2}(\X\times\X)$. Furthermore we have for all $x,y\in\X$,
\begin{equation}\label{boundR}\|\bar R_2(x,y;\cdot)\|_{\barV_2}\leq c\barW_2(x,y).\end{equation}
Let $h:\;\X\times\X\to\rset$ be a symmetric measurable function such that $\bar\mu(|h|)<\infty$. Denote $\theta=\int\int h(x,y)\mu(dx)\mu(dy)$ and define 
\[\bar h_1(x):=\int h(x,z)\mu(dz)-\theta,\;\;\;\bar h_2(x,y)=h(x,y)-\bar h_1(x)-\bar h_1(y)-\theta,\]
\[\bar G_2(x,y):=\int \int\bar R_2(x,y;dz_1,dz_2)\bar h_2(z_1,z_2),\;\;x,y\in\X.\]
We say that $h$ is degenerate when $\bar h_1$ is identically zero. For $x\in\X$, $\delta_x$ denotes the Dirac measure at $x$.
\begin{lemma}\label{lemGfun}
Assume A\ref{A1}.  Suppose that $\bar h_2\in\L_{\barV_2}$. Then $\bar G_2$ is well-defined, $\bar G_2\in \L_{\barW_2}$, and $|\bar G_2|_{\barW_2}\leq c|\bar h_2|_{\barV_2} $ and for all $x,y\in\X$,
\begin{multline}\label{genPoisson}
\bar h_2(x,y)\\
=\int\left(\delta_{x}(dz_1)-P(x,dz_1)\right)\int\left(\delta_{y}(dz_2)-P(y,dz_2)\right)\bar G_2(z_1,z_2).\end{multline}
If in addition $P^s\bar h_{2}\in\L_{\barV_2}$, then $|P^s\bar G_2|_{\barW_2}\leq c|P^s\bar h_2|_{\barV_2} $ for $s\in\{1,2\}$.
\end{lemma}
\begin{proof}
See Section \ref{prooflemGfun}.
\end{proof}

\begin{rem} 
Equation (\ref{genPoisson}) gives a bivariate Poisson's equation which extends the well known univariate Poisson's equation.
\end{rem}

We introduce the function
\begin{multline*}
\Lambda_2(x_1,x_2;y_1,y_2)=\int\left(\delta_{y_1}(dz_1)-P(x_1,dz_1)\right)\int\left(\delta_{y_2}(dz_2)-P(x_2,dz_2)\right)G_2(z_1,z_2)\\
=\bar G_2(y_1,y_2)-P\bar G_2(x_2,y_1)-P\bar G(x_1,y_2)+P^2\bar G_2(x_1,x_2),\;\;\;x_1,x_2,y_1,y_2\in\X.\end{multline*}
Then (\ref{genPoisson}) can be written as $\bar h_2(x,y)=\Lambda_2(x,y,x,y)$. A specially important property of $\Lambda_2$ that we rely on in the sequel is the following. For any $x,y,u,v\in\X$, it is easy to see that
\begin{equation}\label{martproperty}
\int P(x,dy)\Lambda_2(u,x,v,y)=\int P(u,dv)\Lambda_2(u,x,v,y)=0.\end{equation}

Now suppose that we have $\{h_n:\;\X\times\X\to\rset\}$, a family of symmetric measurable functions such that $\bar\mu(|h_n|)<\infty$.We write $\theta_n$, $\bar h_{n,1}$, $\bar h_{n,2}$, $\bar G_{n,2}$,  and $\Lambda_{n,2}$ to denote respectively the quantities $\theta$, $\bar h_{1}$, $\bar h_{2}$, $\bar G_{2}$,  and $\Lambda_{2}$ defined above  with $h=h_n$.

For $1\leq j\leq\ell\leq n$, we introduce the random variables
\[Q_{n,\ell,j}:=\Lambda_{n,2}\left(X_{j-1},X_{\ell-1},X_j,X_\ell\right).\]
For $j< \ell$,  and by the Markov property and (\ref{martproperty}), we have 
\[\PE\left(Q_{n,\ell,j}\vert \F_{\ell-1}\right)=\int P(X_{\ell-1},dz)\Lambda_{n,2}\left(X_{j-1},X_{\ell-1},X_j,z\right)=0,\]
almost surely. This shows that $\{(\sum_{j=1}^{\ell-1}Q_{n,\ell,j},\;\F_{\ell}),\;2\leq \ell\leq n\}$ is a martingale-difference array. We need the following sequences
\begin{multline*}
w_{n,1}(\ell):=\left\{\sum_{j=1}^\ell w_n(\ell,j)+\sum_{j=\ell}^{n}w_n(j,\ell)\right\},\;\;\varpi_{n,1}(\ell):=w_n(\ell,j)-w_n(\ell-1,j)\;\;\\
\varpi_{n,2}(\ell,j):=w_{n}(\ell,j)-w_{n}(\ell,j-1),\;\;\\
\mbox{ and }\;\;\varpi_{n,3}(\ell,j):=w_n(\ell,j)+w_n(\ell-1,j-1)-w_n(\ell,j-1)-w_n(\ell-1,j).\end{multline*}

\begin{lemma}\label{lemHoef}Assume A\ref{A1} and suppose that  $\bar h_{n,2}\in\L_{\barV_2}$ for each $n\geq 1$.
Then
\begin{multline}\label{Hoeffdecomp}
U_n(h_n)=U_{n,0} + \sum_{\ell=1}^n\left\{w_{n,1}(\ell)\bar h_{n,1}(X_\ell) + w_n(\ell,\ell)Q_{n,\ell,\ell}\right\} \\
+\sum_{\ell=1}^n\sum_{j=1}^{\ell-1}w_n(\ell,j)Q_{n,\ell,j} +\zeta_n,\end{multline}
where $U_{n,0}=\theta_n\sum_{\ell=1}^n\sum_{j=1}^\ell w_n(\ell,j)$, and
\begin{multline*}\zeta_n= \sum_{\ell=1}^n\sum_{j=1}^\ell\varpi_n^{(1)}(\ell,j) \left(P\bar G_{n,2}(X_{\ell-1},X_j)-P^2\bar G_{n,2}(X_{\ell-1},X_{j-1})\right) \\
+ \sum_{\ell=1}^n\sum_{j=1}^\ell\varpi_n^{(2)}(\ell,j)\left(P\bar G_{n,2}(X_{j-1},X_{\ell})-P^2\bar G_{n,2}(X_{\ell-1},X_{j-1})\right) \\
+\sum_{\ell=1}^n\sum_{j=1}^\ell\varpi_n^{(3)}(\ell,j)P^2\bar G_{n,2}(X_{\ell-1},X_{j-1})+ \epsilon_n,\end{multline*}
where
\begin{multline*}
\epsilon_n=\sum_{\ell=1}^n\left\{w_n(\ell,0)P\bar G_{n,2}(X_0,X_\ell)-w_n(\ell-1,0)P^2\bar G_{n,2}(X_0,X_{\ell-1})\right) \\
+ \sum_{j=1}^nw_n(n,j)\left(P^2\bar G_{n,2}(X_n,X_j)-P \bar G_{n,2}(X_n,X_j)\right)\\
+\sum_{\ell=1}^n \left(w_n(\ell-1,\ell)P\bar G_{n,2}(X_{\ell-1},X_\ell)-w_n(\ell,\ell)P\bar G_{n,2}(X_\ell,X_{\ell})\right).
\end{multline*}
\end{lemma}
\begin{proof}
See Section \ref{prooflemHoef}. 
\end{proof}

\begin{rem}
The usefulness of this decomposition comes from the fact that the remainder $\zeta_n$ involves either single summations or difference sequences of the weights $w_n$. As a result, these remainders are typically negligible compared to the other terms in the decomposition and one can easily study the asymptotic behavior of $U_n(h_n)$ by focusing on the linear term $\sum_{\ell=1}^n\left\{w_{n,1}(\ell)\bar h_{n,1}(X_\ell)+w_n(\ell,\ell)Q_{n,\ell,\ell}\right\}$, and the quadratic martingale  $\sum_{\ell=1}^n\sum_{j=1}^{\ell-1}w_n(\ell,j)Q_{n,\ell,j}$. 
\end{rem}


\section{Application: Asymptotic variance estimation}\label{sec:avar}
In this section, we use the martingale approximation of Lemma \ref{lemHoef} to study the asymptotics of lag-windows estimators of  asymptotic variance in time series and we apply the results to Markov chain Monte Carlo.  Let $h:\;\X\to\rset$ be a measurable function such that $\mu(|h|^2)<\infty$. We assume without any loss of generality that $\mu(h)=0$. We are interested in the estimation of the long-run variance (or the asymptotic variance) of $h$ defined as:
\begin{equation}\label{asympvar}
\sigma^2(h)=\textsf{Var}_\mu(h(X_0))+2\sum_{\ell\geq 1}\textsf{Cov}_\mu\left(h(X_0),h(X_\ell)\right),\end{equation}
which plays a role in time series analysis and in Markov Chain Monte Carlo.  A classical estimator for $\sigma^2(h)$ is the lag-windows estimators defined as
\begin{equation}\label{kernelest}
\Gamma_{n,b}^2(h):=\gamma_{n,0}+2\sum_{k=1}^{n-1}w_b(kc_n^{-1})\gamma_{n,k},\end{equation}
where $\gamma_{n,k}:=n^{-1}\sum_{j=1}^{n-k}\left(h(X_j)-\mu_n(h)\right)\left(h(X_{j+k})-\mu_n(h)\right)$ is the $k$-th order sample autocovariance with $\mu_n(h)=n^{-1}\sum_{j=1}^n h(X_j)$, $w_b$ is a weight function (with a parameter $b$) and  $\{c_n,\;n\geq 1\}$ is an increasing sequence of positive numbers.
We refer the reader to \cite{priestley81} for detailed discussion on lag-windows estimators. We consider weight functions with the following properties.

\begin{description}
\item[Assumption W] For $b>0$, $w_b:\;[0,\infty)\to [0,1]$ is a continuous function with support $[0,b]$, of class $\mathcal{C}^2$ on the interval $(0,b)$, such that $w_b(b)=0$ and $w_b(0)=1$.
\end{description}\medskip

This assumption allows for the use of all commonly employed weighting functions, including the Bartlett and Parzen kernels. When $w_b(x)=w(x/b)$, an equivalent parametrization of $\Gamma_{n,b}^2(h)$ is $\Gamma_{n,b}^2(h)=\gamma_{n,0}+2\sum_{k=1}^{n-1}w(k/c_n)\gamma_{n,k}$,
 with $c_n\leftarrow b c_n$.  We impose the following ergodicity assumption.

\debutA
\item \label{A4}
There exist measurable functions $V_k:\X\to [1,\infty)$ ($k=1,2,3$), $V_1\leq V_2$, $V_2^2\leq V_3$, such that $PV_3(x)<\infty$ for all $x\in\X$, and $P$ satisfies the assumptions $\C(1,V_1,V_2)$ and $\C(1,V_2^2,V_3)$. Furthermore there exists $q>1$ such that
\begin{equation}\label{momentV3}
\sup_{n\geq 0} \PE\left(V_3^q(X_n)\right)<\infty.\end{equation}
\finA
\vspace{0.2cm}

A\ref{A4} implies A\ref{A1} with $\barV_2(x,y)=V_1(x)V_1(y)$ and $\barW_2(x,y)=V_2(x)V_2(y)$. Define the partial sums $S_{n,k}:=\sum_{j=k+1}^{k+n} h(X_j)$, and the weight $w_{n,b}(0)=n^{-1}$ and $w_{n,b}(k)=2n^{-1}w_b(kc_n^{-1})$ for $k> 0$.  We can rewrite $\gamma_{n,k}$ as 
\[\gamma_{n,k}=n^{-1}\sum_{j=1}^{n-k}h(X_j)h(X_{j+k}) + n^{-3}(n-k)S_{n,0}^2-n^{-2}S_{n,0}\left(S_{n-k,0} + S_{n-k,k}\right),\]
 so that 
\begin{equation}\label{kernelest2}
\Gamma_{n,b}^2(h)=\sum_{\ell=1}^n\sum_{j=1}^\ell w_{n,b}(\ell-j)h(X_j)h(X_\ell) + R_n,\;\;
\end{equation}
where
\begin{multline}\label{termRn}
R_n=2n^{-2}S_{n,0}^2\sum_{k=1}^{n-1}w_b(kc_n^{-1}) \left(1-\frac{k}{n}\right) \\
-2n^{-2}S_{n,0}\left(\sum_{j=2}^n h(X_j)\sum_{k=1}^{j-1}w_b(kc_n^{-1}) 
+ \sum_{j=1}^{n-1} h(X_j)\sum_{k=1}^{n-j}w_b(kc_n^{-1})\right) -n^{-2}S_{n,0}^2.\end{multline}
If we set aside the term $R_n$, the expression (\ref{kernelest2}) is of the form (\ref{Un}) with $h_n(x,y)=h(x)h(y)$ and $w_n(\ell,j)=w_{n,b}(\ell-j)$. Here we have $h_{n,1}(x)=\int h(x)h(y)\mu(dy)=0$, $\theta_n=0$, and $h_{n,2}(x,y)=h(x)h(y)$. Define 
\[G(x):=\sum_{j\geq 0} P^j h(x),\;\;\mbox{ and }\;\;PG(x)=\int P(x,dz)G(z),\;\;x\in\X.\]
Then $\bar G_2(x,y)=G(x)G(y)$, $P\bar G_2(x,y)=PG(x) G(y)$, and $P^2\bar G_2(x,y)=PG(x)PG(y)$. Therefore 
\[Q_{n,\ell,j}=Q_\ell Q_j,\;\;\mbox{ where }\;\;Q_\ell=G(X_\ell)-PG(X_{\ell-1}).\]
As above, $\{(Q_\ell,\F_\ell),\;\ell\geq 1\}$ is a martingale: $\PE\left(Q_\ell\vert \F_{\ell-1}\right)=0$. From Lemma \ref{lemHoef} we obtain the following.

\begin{theo}\label{thm3}
Assume (A\ref{A4}) and (W) and $h\in\L_{V_1}$. For all $n\geq 1$,
\begin{equation}\label{decompGam}
\Gamma_{n,b}^2(h)=n^{-1}\sum_{\ell=1}^n Q_\ell^2+\sum_{\ell=1}^n\sum_{j=1}^{\ell-1} w_{n,b}(\ell-j)Q_\ell Q_j +R_n + \zeta_n.\end{equation}
Furthermore, there exist $p>1$ and a finite constant $c$ such that for all $n\geq 3$,
\begin{multline*}
\PE^{1/p}\left(|\zeta_n|^p\right)\leq c c_n^{-1+\frac{1}{2}\vee\frac{1}{p}},\;\;\;\PE^{1/p}\left(\left|R_n\right|^p\right)\leq c n^{-1}c_n,\\
\mbox{ and }\;\;\;\PE^{1/p}\left(\left|\sum_{\ell=1}^n\sum_{j=1}^{\ell-1} w_{n,b}(\ell-j)Q_\ell Q_j\right|^p\right)\leq c \left(\frac{c_n}{n}\right)^{\frac{1}{2}}n^{-\frac{1}{2}+\frac{1}{p}\vee\frac{1}{2}}.\end{multline*}
\end{theo}
\begin{proof}
See Section \ref{proofthm3}.
\end{proof}


A clearer picture of the behavior of the lag-window estimator  emerges from this result. For $p\geq 2$, we have
\begin{equation}\label{FulldecompGam}
\Gamma_{n,b}^2(h) = \underset{O_p(1)}{\underbrace{n^{-1}\sum_{\ell=1}^n Q_\ell^2}}
 + \underset{O_p\left(\sqrt{\frac{c_n}{n}}+\frac{c_n}{n}\right)}{\underbrace{\sum_{\ell=1}^n\sum_{j=1}^{\ell-1} w_{n,b}(\ell-j)Q_\ell Q_j + R_n}}
 + \underset{O_p(c_n^{-1/2})}{\underbrace{\zeta_n}},
\end{equation}
By the law of large numbers for Markov chain the term $n^{-1}\sum_{\ell=1}^n Q_\ell^2$ converges to $\sigma^2(h)$.  As the result, Theorem \ref{thm3} implies that $\Gamma_{n,b}^2(h)$ converges in probability to $\sigma^2(h)$ provided $c_n\to\infty$, $c_n=o(n)$ and $p\geq 2$ (for $1<p<2$, specific rate assumption on $c_n$ might be needed). The decomposition (\ref{FulldecompGam}) also gives some insight into the well known fact that $\Gamma_{n,b}(h)$ often has poor finite-sample properties in estimating $\sigma^2(h)$, particularly for highly correlated time-series.  Indeed, for $c_n=o(n)$, both terms $R_n + \sum_{\ell=1}^n\sum_{j=1}^{\ell-1} w_{n,b}(\ell-j)Q_\ell Q_j$ and $\zeta_n$ converge to zero but at antagonistic rates. If $c_n\approx n$, then $\zeta_n\approx O_P(n^{-1/2})$ but then  $R_n+\sum_{\ell=1}^n\sum_{j=1}^{\ell-1} w_{n,b}(\ell-j)Q_\ell Q_j\approx O(1)$. Whereas for $c_n\ll n$, the convergence of $\zeta_n$ is slow ($\zeta_n=O_P(c_n^{-1/2})$) but $R_n+\sum_{\ell=1}^n\sum_{j=1}^{\ell-1} w_{n,b}(\ell-j)Q_\ell Q_j$ vanishes quickly. 

When the goal is to construct confidence interval for $\mu(h)$ (and one is not interested in estimating $\sigma^2(h)$ per se), it has been suggested to use the lag-window estimator $\Gamma_{n,b}^2(h)$ with $c_n=n$, the so-called ``fixed-b asymptotics'' (\cite{neave70, kieferetvogelsang05}). With $c_n=n$, $\Gamma_{n,b}^2(h)$ no longer converges to $\sigma^2(h)$, but as it turns out, asymptotically valid confidence intervals can still be derived for $\mu(h)$. We have the following.

%
%

\begin{theo}\label{thm4}
Under the assumption of Theorem \ref{thm3}, the following holds true.
\begin{enumerate}
\item If $p\geq 2$ and $c_n=o(n)$, then $\Gamma_{n,b}^2(h)$ converges in probability to $\sigma^2(h)$. Furthermore, assuming $\Gamma_{n,b}^2(h)> 0$ almost surely, 
\[\left\{n\Gamma_{n,b}^2(h)\right\}^{-1/2}\sum_{j=1}^n \left(h(X_j)-\mu(h)\right)\stackrel{w}{\to} \mathcal{N}(0,1).\]
\item Let $\{B(t),\,0\leq t\leq 1\}$ be the standard Browian motion. If $c_n=n$, then $\Gamma_{n,b}^2(h)\stackrel{w}{\to} \sigma^2(h)\textsf{K}_b$, where
\begin{multline*}
\textsf{K}_b=1+2\int_0^1\int_0^tw_b(t-s)dB(s)dB(t) \\
- 2B(1)\int_0^1g_b(t)dB(t)+2B^2(1)\int_0^1(1-t)w_b(t)dt,\end{multline*}
where $g_b(t)=\int_0^tw_b(u)du +\int_0^{1-t}w_b(u)du$. Furthermore, assuming $\Gamma_{n,b}^2(h)>0$ almost surely,
\[\left\{n \Gamma_{n,b}^2(h)\right\}^{-1/2}\sum_{j=1}^n \left(h(X_j)-\mu(h)\right)\stackrel{w}{\to}\frac{B(1)}{\sqrt{\textsf{K}_b}}.\]
\end{enumerate}
\end{theo}
\begin{proof}
See Section \ref{proofthm4}.
\end{proof}

\vspace{0.2cm}
By Theorem \ref{thm4}~(1) an asymptotically valid $(1-\alpha)$-confidence interval for $\mu(h)$ is 
\begin{equation}\label{CI1}
\mu_n(h)\pm z_{1-\alpha/2}\frac{\hat\sigma_n(h)}{\sqrt{n}},\end{equation}
where $z_{1-\alpha/2}$ is the $(1-\alpha/2)$-quantile of the standard normal distribution and where $\hat\sigma_n(h)=\sqrt{\Gamma_{n,b}^2(h)}$, with $b=1$, $c_n=o(n)$. Typical choice of $c_n$ includes $c_n=n^{-\delta}$, $\delta\in (0,1)$ typically around $0.5$. Theorem \ref{thm3} (2) provides another asymptotically valid confidence interval for $\mu(h)$:
\begin{equation}\label{CI2}
\mu_n(h)\pm t_{1-\alpha/2}\frac{\tilde\sigma_n(h)}{\sqrt{n}},\end{equation}
where $t_{1-\alpha/2}$ is the $(1-\alpha/2)$-quantile of the distribution of $B(1)/\sqrt{\textsf{K}_b}$ and where $\tilde\sigma_n(h)=\sqrt{\Gamma_{n,b}^2(h)}$, with $c_n=bn$, with $b\in (0,1)$. 

Although the limiting distribution $B(1)/\sqrt{\textsf{K}_b}$ is non-standard, it can be simulated, for example by Euler discretization of the stochastic integrals in $\textsf{K}_b$. We report  in Table 1 the $95\%$ quantiles of the distribution of  $B(1)/\sqrt{\textsf{K}_b}$ using $w_b(x)=\textbf{1}_{(0,b)}(x)$, $w_b(x)=(1-x/b)\textbf{1}_{(0,b)}(x)$ and $w_b(x)=(1-(x/b)^2)\textbf{1}_{(0,b)}(x)$, and for different values of $b$, based on $10,000$ replications of $B(1)/\sqrt{\textsf{K}_b}$. The distribution departs further from the standard normal distribution as $b$ increases.

\begin{table}[h]\label{table1}
\begin{center}
\small
\begin{tabular}{cccc}
\hline
&$w_b(x)=1-x/b$& $w_b(x)=1-(x/b)^2$&$w_b(x)=\textbf{1}_{(0,1)}(x/b)$\\
\hline
$b=0.3$ & $2.828$ & $4.134$ & $5.496$ \\
$b=0.5$ & $3.557$ & $6.580$ & $6.299$ \\
$b=0.9$ & $4.735$ & $12.575$ & $13.045$ \\
\hline
\end{tabular}
\caption{$0.975$-quantile of the distribution of $B(1)/\sqrt{\textsf{K}_b}$.}
\end{center}
\end{table}

In the next simulation examples, we compare the finite sample properties of these two confidence intervals in terms of coverage probability and interval length. All the simulations are performed using the Bartlett kernel $w(x)=1-x$.

%
%

\subsection{Illustration: the Garch$(1,1)$ model}
Consider the linear GARCH$(1,1)$ model defined as follows. $h_0\in(0,\infty)$, $u_0\sim \mathcal{N}(0,h_0)$ and for $n\geq 1$
\begin{eqnarray*}
u_n&=&h_n^{1/2}\epsilon_n\\
h_n&=&\omega +\beta h_{n-1} + \alpha u_{n-1}^2,\end{eqnarray*}
where $\{\epsilon_n,\;n\geq 0\}$ is i.i.d. $\mathcal{N}(0,1)$ and $\omega>0$, $\alpha\geq 0$, $\beta\geq 0$.  We assume that $\alpha,\beta$ satisfy

\debutC
\item There exists $\nu>0$ such that
\begin{equation}\label{E1}
\E\left[\left(\beta+\alpha Z^2\right)^\nu\right]<1,\;\;\;Z\sim\mathcal{N}(0,1).\end{equation}
\finC

It is shown by \cite{meitzetsaikkonen08} (Theorem 2) that under (\ref{E1}) the joint process $\{(u_n,h_n),\;n\geq 0\}$ is a phi-irreducible aperiodic Markov chain that admits an invariant distribution  and is geometrically ergodic with a drift function $V(u,h)=1+h^\nu+|u|^{2\nu}$. Therefore for $\nu\geq 2$, A\ref{A4} holds with $V_1=V_2=V^{1/2}$, and $V_3=V$. We are interested in a confidence interval for $\mu(h)$ where  $h(u)=u^2$ which belongs to $\L_{V_1}$. The exact value is $\mu(h)=\omega(1-\alpha-\beta)^{-1}$.

For the simulations we set $\omega=1,\;\alpha=0.1,\;\beta=0.7$ which gives $\mu(h)=5$. We compare the confidence intervals (\ref{CI1}) and (\ref{CI2}) by computing (by Monte Carlo) their coverage probabilities and average lengths. The comparison is performed using sample paths of length $60,000$ from the GARCH$(1,1)$ Markov chain.  The results are plotted in Figure 1 and shows across the board better coverage probability of the fixed-b confidence interval but, as expected,  at the expense of a slightly wider confidence intervals. 

\vspace{0.1cm}
\begin{center}
\scalebox{0.4}{\includegraphics{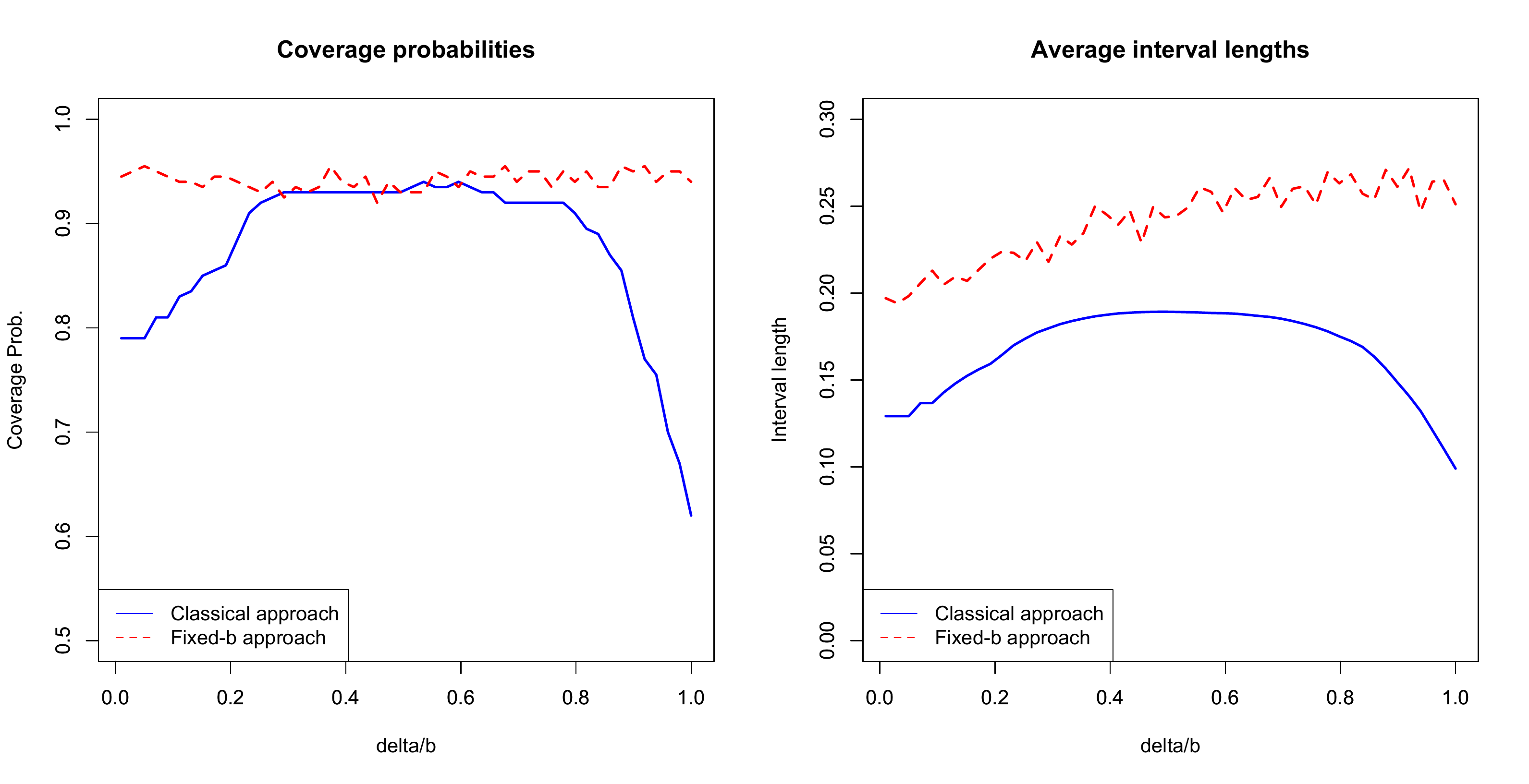}}\\
\noindent\underline{Figure 1}: Coverage probabilities plots for various values of $\delta$ (classical confidence interval) and $b$ (fixed-b confidence interval).
\end{center}

\subsection{Markov Chain Monte Carlo}
Markov Chain Monte Carlo (MCMC) is a popular computational tools to obtain random samples from intractable and high-dimensional distributions (see e.g. \cite{robertsetrosenthalsurvey} for a survey and for additional references). 

Suppose that we interested in sampling from the probability measure $\mu$ and compute the integral $\mu(h)=\int h(x)\mu(dx)$. Let $\{X_n,\;n\geq0\}$ be a Markov chains with transition kernel $P$, invariant distribution $\mu$ and initial distribution $\rho$. By simulating the Markov chain, we approximate $\mu(h)$ by the Monte Carlo average $\mu_n(h)=n^{-1}\sum_{k=1}^n h(X_k)$. Furthermore, under A\ref{A4}, $\lim_{n\to\infty} n^{1/2} \textsf{Var}\left(\mu_n(h)\right)=\sigma^2(h)$, as given by (\ref{asympvar}), and a central limit theorem holds: $n^{-1/2}\sum_{k=1}^n\left(h(X_k)-\pi(h)\right)\stackrel{w}{\to} N(0,\sigma^2(h))$.  Therefore (\ref{CI1}) and (\ref{CI2}) provide two valid confidence intervals for $\mu(h)$. We compare the coverage probabilities and average interval lengths of these two confidence interval procedures with the following simulation example.

\subsubsection{Illustration: a Poisson regression model}
We undertake the comparison using a log-linear model taken from \cite{gelmanetal03}. For $e=1,\ldots,N_e$ and $p=1,\ldots,N_p$, the variables $y_{ep}$ are conditionally independent given $(\{\beta_p\}$, $\{\eps_{ep}\})\in\rset^{N_p}\times\rset^{N_eN_p}$, with conditional distribution 
\begin{equation}y_{ep}\sim \mathcal{P}\left(n_{ep}e^{\mu+\alpha_e+\beta_p+\eps_{ep}}\right),\;\; e=1,\ldots,N_e,\;\; p=1,\ldots,N_p,\end{equation}
where $\mathcal{P}(\lambda)$ is the Poisson distribution with
parameter $\lambda$. In the above display, $\{n_{ep}\}$ is a deterministic baseline covariate, and $\mu\in\rset$, $\{\alpha_e\}\in\rset^{N_e}$ are parameters. We assume that $\{\beta_p\}$ and $\{\eps_{ep}\}$ are independent with distributions
\begin{equation}\beta_p\stackrel{iid}{\sim}N(0,\sigma^2_\beta),\;\;\;\;\eps_{ep}\stackrel{iid}{\sim}N(0,\sigma^2_\eps),\;\; e=1,\ldots,N_e,\;\; p=1,\ldots,N_p,\end{equation}
for some parameters $\sigma_\beta^2>0$, $\sigma_\eps^2>0$. We assume a diffuse prior for
$(\mu,\alpha,\sigma_\beta^2,\sigma_\eps^2)$ ($\sigma_\epsilon^2>0,
\sigma_\beta^2>0$) with the additional constraint that
$\alpha_{N_e}=-\sum_{k=1}^{N_e-1}\alpha_k$. Let
$\theta=(\mu,\alpha,\beta,\epsilon,\sigma^2_\epsilon,\sigma^2_\beta)\in\rset^{3+N_e-1+(N_p+1)N_e}$. The
posterior distribution of $\theta$ given $\mathcal{D}=(y_{ep},n_{ep})$
takes the form
\begin{multline}\label{postEx}
\pi(\theta\vert \mathcal{D})\propto \exp\left\{\sum_{e,p}y_{e,p}(\mu+\alpha_e+\beta_p+\epsilon_{e,p})-n_{ep}e^{\mu+\alpha_e+\beta_p+\epsilon_{ep}}\right.\\
\left.-\frac{N_eN_p}{2}\log\sigma_\epsilon^2-\frac{N_p}{2}\log\sigma^2_\beta-\frac{1}{2\sigma_\epsilon^2}\sum_{e,p}\epsilon^2_{e,p}
-\frac{1}{2\sigma^2_\beta}\sum_{p=1}^{N_p}\beta_p^2\right\}.\end{multline}
This posterior distribution is typical of probability distributions for which MCMC is useful. We set $N_e=3$ and $N_p=20$. Suppose that we are interested in a confidence interval for the posterior mean of the parameter $\alpha_1$, i.e. $\int \alpha_1\pi(\theta\vert \mathcal{D})d\theta$. To compare the two confidence intervals methods described above, we generate an artificial dataset with $(\alpha_1,\alpha_2,\mu,\sigma_\eps^2,\sigma_\beta^2)=(0.35, 0.15, -1.0, 0.1, 0.3)$. We run a  preliminary MCMC sampler for 6 millions ($6\times 10^6$) iterations and compute its sample mean. We obtain $\bar \alpha_1=0.3309$. We take this value to be $\int \alpha_1\pi(\theta)d\theta$.

To compare the two confidence interval methods, we use a Random Walk Metropolis (RWM) algorithm with proposal kernel $\mathcal{N}(0,\kappa\Sigma)$ where  $\kappa$ and $\Sigma$ are selected (from a preliminary simulation) to yield a reasonably good mixing of the chain. We run the MCMC sampler for $60,000$ iterations and discard the first $10,000$ iterations as burn-in. We repeat the simulations $200$ times in order to estimate the coverage probabilities and interval lengths. The results are given in Figure 2. We find again that in terms of finite sample behavior, the fixed-b confidence interval is more robust to the choice of $b$.

\vspace{0.1cm}
\begin{center}
\scalebox{0.4}{\includegraphics{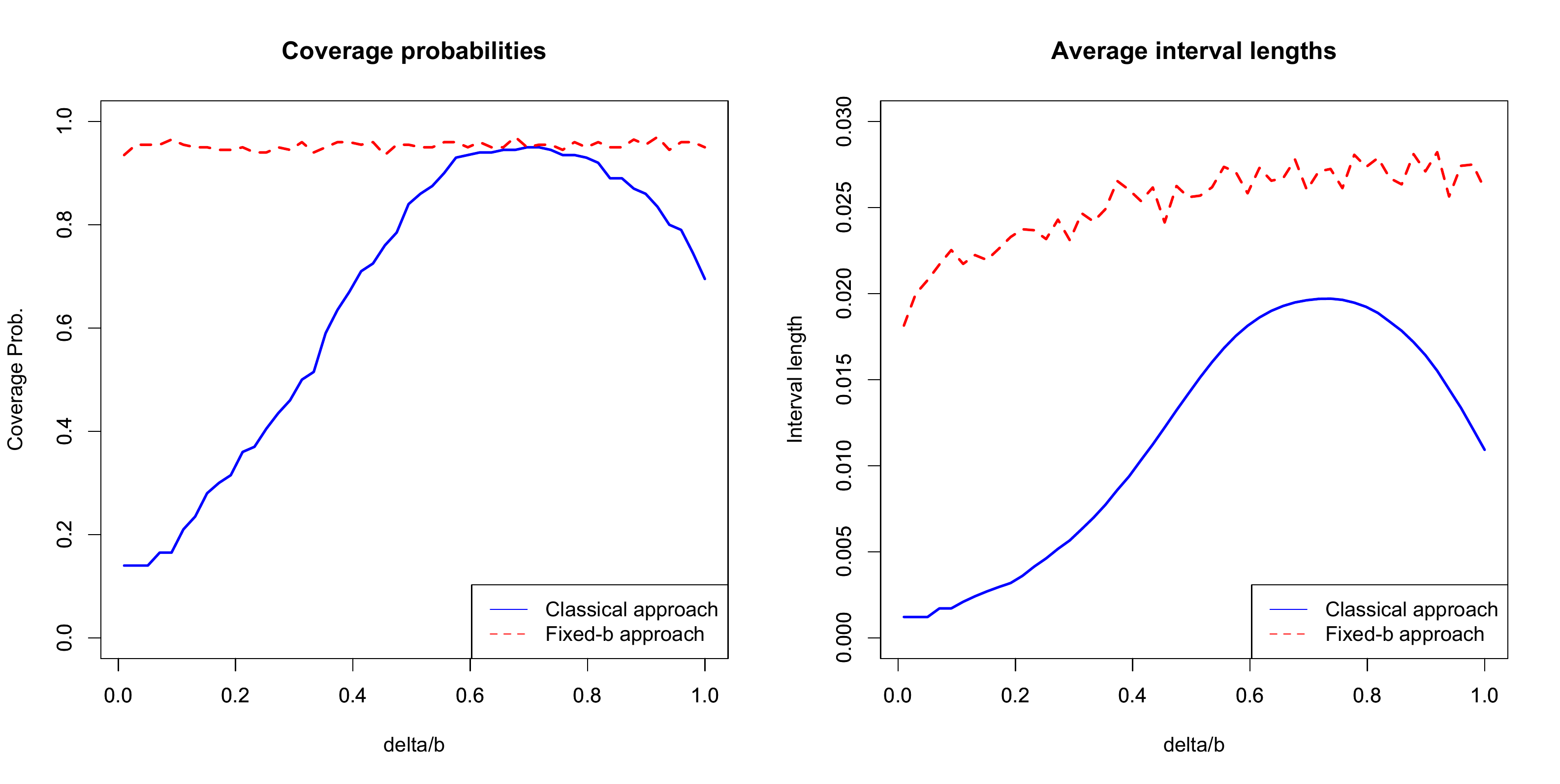}}\\
\noindent\underline{Figure 2}: Coverage probability and confidence interval length for $\alpha_1$ and for different values of $b$ and $\delta$.
\end{center}

%

\section{Application: a CLT for U-statistics with varying kernels}\label{sec:Ustat}
 U-statistics with varying kernels are a special case of quadratic forms and correspond to setting $w_n(\ell,\ell)=0$ and $w_n(\ell,j)=1$ if $\ell\neq j$. We thus have  \[U_n(h_n)=\sum_{\ell=2}^n\sum_{j=1}^{\ell-1}h_n(X_\ell,X_j).\]
 We illustrate another application of Lemma \ref{lemHoef}, by deriving a CLT for $U_n(h_n)$. U-statistics and U-statistics with varying kernels play an important role in nonparametric and semi-parametric statistics.   In the present case,
 under A\ref{A1}, Lemma \ref{lemHoef} reduces to
 \[U_n(h_n)={n\choose 2}\theta_n +(n-1)\sum_{\ell=1}^n \bar h_{n,1}(X_\ell) +\sum_{\ell=2}^n\sum_{j=1}^{\ell-1}Q_{n,\ell,j} +\zeta_n.\]

We impose the following moment assumption
\debutB
\item \label{B1} With $\barV_2$ and $\barW_2$ as in A\ref{A1}, suppopse that $P\barV_2\leq c\barV_2$, $P\barW_2^2\leq c\barW_2^2$, 
 and for $p=2$,
\begin{equation}\label{momentB1}\sup_{\ell,j\geq 0}\PE\left(\barW_2^p(X_\ell,X_j)\right)<\infty.\end{equation}
 \finB

We recall from Lemma \ref{lemGfun} that  for $s\in\{1,2\}$, $|P^s\bar G_{n,2}(x,y)|\leq c|P^s\bar h_{n,2}|_{\barV_2}\barW_2(x,y)$, and 
\begin{multline*}|P^2\bar h_{n,2}|_{\barV_2}=\sup_{x,y\in\X}\{\barV_2(x,y)\}^{-1}\left|\int P(x,du)\int P(y,dv)\bar h_{n,2}(u,v)\right|\\
\leq |P\bar h_{n,2}|_{\barV_2}\sup_{x,y\in\X}\{\barV_2(x,y)\}^{-1}\int P(x,du)\barV_2(u,y)\leq c|P\bar h_{n,2}|_{\barV_2}\leq c\nnorm{\bar h_{n,2}}_{p,\barV_2},\end{multline*}
for all $p\geq 1$, using the assumption $P\barV_2\leq c\barV_2$ in B\ref{B1}. In combination with (\ref{momentB1}), and the expression of $\zeta_n$ in Lemma \ref{lemHoef}, it follows that
\begin{equation}\label{boundzeta}
\PE^{1/2}\left(|\zeta_n|^2\right)\leq c\,n\left( |P\bar h_{n,2}|_{\barV_2} + |P^2\bar h_{n,2}|_{\barV_2}\right)\leq c\,n\,\nnorm{\bar h_{n,2}}_{2,\barV_2}.\end{equation}
 By definition, $Q_{n,\ell,j}=\bar h_{n,2}(X_\ell,X_j)-P\bar G_{n,2}(X_{\ell-1},X_j)-P\bar G_{n,2}(X_{j-1},X_\ell)+P^2 \bar G_{n,2}(X_{\ell-1},X_{j-1})$ and assuming that $P\bar h_{n,2},\in\L_{\barV_2}$, one obtains
\begin{multline}\label{boundQ}
|Q_{n,\ell,j}|\leq \left|\bar h_{n,2}(X_\ell,X_j) \right|+ c|P \bar h_{n,2}|_{\barV_2}\left(\barW_2(X_{\ell-1},X_j)+\barW_2(X_{j-1},X_\ell)\right)\\
+ c|P^2 \bar h_{n,2}|_{\barV_2}\barW_2(X_{\ell-1},X_{j-1}).\end{multline}
Thus for $j<\ell$,
\begin{multline*}
\PE\left(Q_{n,\ell,j}^2\vert \F_{\ell-1}\right)\leq 4\int P(X_{\ell-1},dz)\left|\bar h_{n,2}(z,X_j) \right|^2 \\
+ 4\nnorm{\bar h_{n,2}}_{2,\barV_2}\PE\left(\barW_2^2(X_{\ell-1},X_j) +\barW_2^2(X_{\ell},X_{j-1})+ \barW_2^2(X_{\ell-1},X_{j-1})\vert \F_{\ell-1}\right).\end{multline*}
Taking the expectation on both side and using (\ref{momentB1}), it follows that for all $n\geq 1$,
\begin{equation}\label{momentQ}
\PE^{1/2}\left(Q_{n,\ell,j}^2\right)\leq c\nnorm{\bar h_{n,2}}_{2,\barV_2}.\end{equation}

 We impose an addition stability assumptions.

\debutB
\item \label{B2} With $\barW_2$ as in A\ref{A1}, there exists measurable functions $\overline{\mathcal{U}}_1\leq \overline{\mathcal{V}}_1:\,\X\to [1,\infty)$, a symmetric measurable function $\overline{\mathcal{U}}_2:\,\X\times\X\to [1,\infty)$ such that $P$ satisfies $\C(1,\overline{\mathcal{U}}_1,\overline{\mathcal{V}}_1)$ and for all $m\geq 0$, all $x_1,x_2,x_3\in\X$,
\begin{multline}\label{controlA2}
\int P(x_1,dz_1) \int P^{m}(z_1,dz_2)\\
\times\left(\barW_2(z_2,z_1)+\barW_2(z_2,x_1)\right)\left(\barW_2(z_2,x_2)+\barW_2(z_2,x_3)\right)\\
\leq c\,\overline{\mathcal{U}}_1(x_1)\overline{\mathcal{U}}_2(x_2,x_3).\end{multline}
Futhermore, 
\[\sup_{\ell\geq 1} \PE\left(\overline{\mathcal{V}}_1(X_\ell)\overline{\mathcal{U}}(X_\ell,X_{\ell-1})\right)<\infty.\]
\finB

\vspace{0.2cm}
\begin{rem}
Assumptions B\ref{B1}-B\ref{B2} are similar to the assumptions imposed in \cite{dehlingetal10} (Theorem~1.8) to obtain a CLT for U-statistics of stationary dependent processes. Assumption B\ref{B2} can be easy to check. For example if $\barW_2$ is given by $\barW_2(x,y)=W(x)W(y)$ and $P^mW\leq cW$ for all $m\geq 0$, then (\ref{controlA2}) holds with $\overline{\mathcal{U}}_1(x)=PW^2(x)+W(x)PW(x)$ and  $\overline{\mathcal{U}}_2(x,y)=W(x)+W(y)$.
\end{rem}
Define \begin{multline*}
\sigma_{n,1}^2:=\int\{\mu P\}(dx,dy)L_n^2(x,y)\\
=\textsf{Var}_\mu\left(h_{n,1}(X_0)\right)+2\sum_{\ell\geq 1}\textsf{Cov}_\mu\left(h_{n,1}(X_0),h_{n,1}(X_\ell)\right),\;\;\mbox{ and }\;\;
\sigma_n^2\,:=\,n(n-1)^2\sigma_{n,1}^2.\end{multline*}

\begin{theo}\label{theoclt1}
Assume A\ref{A1}, B\ref{B1}-B\ref{B2} and let $\{h_n,\,n\geq 1\}$ be such that $\bar h_{n,2}\in\L_{\barV_2}$. Suppose also that 
\begin{equation}\label{condeq1}
\nnorm{\bar h_{n,2}}_{2,\barV_2}=o\left(n^{1/2}\sigma_{n,1}\right).\end{equation}
Then
\[\frac{1}{\sigma_{n}}\left(U_n(h_n)-\theta_n{n\choose 2} \right)-\frac{1}{\sigma_{n,1}\sqrt{n}}\sum_{\ell=1}^n\bar h_{n,1}(X_\ell)\to 0,\;\;\mbox{ in probab. }\]
\end{theo}
\begin{proof}
See Section \ref{prooftheoclt1}.
\end{proof}

\begin{rem}
From the above result, it is clear that if $\frac{1}{\sigma_{n,1}\sqrt{n}}\sum_{\ell=1}^n\bar h_{n,1}(X_\ell)$ converges weakly to $\mathcal{N}(0,1)$, then so does $\frac{1}{\sigma_{n}}\left(U_n(h_n)-\theta_n{n\choose 2} \right)$. If $h_n$ does not depend on $n$, (\ref{condeq1}) automatically holds and Theorem \ref{theoclt1} implies  a standard CLT for U-statistics (\cite{yoshihara76,dehlingetal10}). But unlike these previous works, Theorem \ref{theoclt1} does not assume stationarity and for Markov chains, the weak dependence assumption A\ref{A1} is sometimes easier to check than mixing assumptions.

The theorem describes the limiting behavior of $U_n(h_n)$ in the case where the kernels $h_n$ are not degenerate and the quadratic term $\sum_{\ell=1}^n\sum_{j=1}^{\ell-1} Q_{n,\ell,j}$ is negligible. In general, the quadratic term needs not be negligible. In which case a  correct account of the limiting behavior of $U_n(h_n)$ will then require a  joint study the processes $\sum_{\ell=1}^n\bar h_{n,1}(X_\ell)$ and $\sum_{\ell=1}^n\sum_{j=1}^{\ell-1}Q_{n,\ell,j}$. 
\end{rem}

\section{Proofs}\label{sec:proofs}
\subsection{Proof Lemma \ref{lemGfun}}\label{prooflemGfun}

That $\bar G_2\in \L_{\barW_2}$, and $|\bar G_2|_{\barW_2}\leq c|\bar h_2|_{\barV_2} $ follows from (\ref{boundR}). Set 
\[\pi_{n,m}(x,y;(du,dv))=\left(P^n(x,du)-\mu(du)\right)\bigotimes\left(P^m(y,dv)-\mu(dv)\right).\]
 Since $P^s\barW_2(x,y)<\infty$ for all $x,y\in\X$ and $s\in\{1,2\}$ by A\ref{A1}, we deduce that the rhs of (\ref{genPoisson}) is well-defined and can be written as $ \bar G_2(x,y)-P\bar G_2(y,x)-P\bar G_2(x,y)+P^2\bar G_2(x,y)$. By dominated convergence,
 \begin{multline*}
 \bar G_2(x,y)-P\bar G_2(y,x)-P\bar G_2(x,y)+P^2\bar G_2(x,y)\\
 =\lim_{N,M\to\infty}\sum_{n=0}^N\sum_{m=0}^M\int\left(\delta_{x}(dz_1)-P(x,dz_1)\right)\int\left(\delta_{y}(dz_2)-P(y,dz_2)\right)\pi_{n,m}\bar h_2(z_1,z_2),\\
 =\lim_{N,M\to\infty}\left\{\bar h_2(x,y)-\pi_{N+1,0}\bar h_2(x,y)-\pi_{0,M+1}\bar h_2(x,y)+\pi_{N+1,M+1}\bar h_2(x,y)\right\}\\
 =\bar h_2(x,y),\end{multline*}
proving (\ref{genPoisson}). The bound $|P^s\bar G_2|_{\barW_2}\leq c|P^s\bar h_2|_{\barV_1} $ is obtained by showing in a similar way that
\begin{multline*}
P^s\bar G_2(x,y)=\lim_{N,M\to\infty}\sum_{n=0}^N\sum_{m=0}^M \pi_{n,m}\{P^s\bar h_2\}(x,y)\\
=\int\int\bar R_2(x,y;dz_1,dz_2)\{P^s\bar h_2\}(z_1,z_2).\end{multline*}

\begin{flushright}
$\square$
\end{flushright}

\subsection{Proof Lemma \ref{lemHoef}}\label{prooflemHoef}

From the definition, we have $h_n(x,y)=\theta_n+\bar h_{n,1}(x)+\bar h_{n,1}(y)+\bar h_{n,2}(x,y) $, and we deduce after some rearrangements that
\[U_n(h_n)= U_{n,0} + \sum_{\ell=1}^nw_{n,1}(\ell)\bar h_{n,1}(X_\ell) +\sum_{\ell=1}^n\sum_{j=1}^\ell w_{n}(\ell,j)\bar h_{n,2}(X_\ell,X_j).\]
Using (\ref{genPoisson}), we write 
\begin{multline*}
\bar h_{n,2}(X_\ell,X_j)= \Lambda_{n,2}(X_j,X_\ell,X_j,X_\ell)\\
=Q_{n,\ell,j} + \Lambda_{n,2}(X_j,X_\ell,X_j,X_\ell)-\Lambda_{n,2}(X_{j-1},X_{\ell-1},X_j,X_\ell)\\
=Q_{n,\ell,j} + \left(P\bar G_{n,2}(X_{\ell-1},X_j)-P\bar G_{n,2}(X_\ell,X_j)\right) \\
+ \left(P\bar G_{n,2}(X_{j-1},X_\ell)-P\bar G_{n,2}(X_j,X_\ell)\right)  + \left(P^2\bar G_{n,2}(X_j,X_{\ell})-P^2\bar G_{n,2}(X_{j-1},X_{\ell-1})\right).\end{multline*}
Rearranging the terms, it is easy to verify that
\begin{multline*}\label{proofdecompeq1}
\sum_{\ell=1}^n\sum_{j=1}^\ell w_n(\ell,j)\bar h_{n,2}(X_\ell,X_j)=\sum_{\ell=1}^n\sum_{j=1}^\ell w_n(\ell,j)Q_{n,\ell,j} \\
+ \sum_{\ell=1}^n\sum_{j=1}^\ell\varpi_n^{(1)}(\ell,j) \left(P\bar G_{n,2}(X_{\ell-1},X_j)-P^2\bar G_{n,2}(X_{j-1},X_{\ell-1})\right) \\
+ \sum_{\ell=1}^n\sum_{j=1}^\ell\varpi_n^{(2)}(\ell,j)\left(P\bar G_{n,2}(X_{j-1},X_\ell)-P^2\bar G_{n,2}(X_{j-1},X_{\ell-1})\right) \\
+\sum_{\ell=1}^n\sum_{j=1}^\ell\varpi_n^{(3)}(\ell,j)P^2\bar G_{n,2}(X_{j-1},X_{\ell-1})+ \epsilon_n,\end{multline*}
where $\epsilon_n$ is comprised of the remainder telescoping sums. We obtain
\begin{multline*}
\epsilon_n=\sum_{\ell=1}^n\left\{w_n(\ell,0)P\bar G_{n,2}(X_0,X_\ell)-w_n(\ell-1,0)P^2\bar G_{n,2}(X_0,X_{\ell-1})\right) \\
+ \sum_{j=1}^nw_n(n,j)\left(P^2\bar G_{n,2}(X_n,X_j)-P \bar G_{n,2}(X_n,X_j)\right)\\
+\sum_{\ell=1}^n \left(w_n(\ell-1,\ell)P\bar G_{n,2}(X_{\ell-1},X_\ell)-w_n(\ell,\ell)P\bar G_{n,2}(X_\ell,X_{\ell})\right).
\end{multline*}

\begin{flushright}
$\square$
\end{flushright}

\subsection{Proof theorem \ref{thm3}}\label{proofthm3}
A\ref{A4} implies that we can find $p>1$ such that
\begin{equation}\label{momentV4}\sup_{n\geq 0}\PE\left(V_2^{2p}(X_n)\right)<\infty.\end{equation}
From (\ref{kernelest2}), $\Gamma_{n,b}^2(h)=\sum_{\ell=1}^n\sum_{j=1}^\ell w_{n,b}(\ell-j)h(X_j)h(X_\ell) + R_n$, where $w_{n,b}(0)=n^{-1}$ and $w_{n,b}(\ell)=2n^{-1}w_b(kc_n^{-1})$ (in particular $w_{n,b}(\ell)=0$ for $\ell<0$). Set $\Delta_{n,b}^{(1)}(\ell)=w_{n,b}(\ell)-w_{n,b}(\ell-1)$ and $\Delta_{n,b}^{(2)}(\ell)=2w_{n,b}(\ell)-w_{n,b}(\ell+1)-w_{n,b}(\ell-1)$. Then  Lemma \ref{lemHoef} applied to $\sum_{\ell=1}^n\sum_{j=1}^\ell w_{n,b}(\ell-j)h(X_j)h(X_\ell)$ gives:
\[\Gamma_{n,b}^2(h)=n^{-1}\sum_{\ell=1}^n Q_\ell^2+\sum_{\ell=1}^n\sum_{j=1}^{\ell-1} w_{n,b}(\ell-j)Q_\ell Q_j +R_n + \zeta_n,\]
where
\begin{multline*}
\zeta_n=\sum_{\ell=1}^n PG(X_{\ell-1})\sum_{j=1}^{\ell}\Delta_{n,b}^{(1)}(\ell-j)Q_j + \sum_{\ell=1}^n Q_\ell \sum_{j=1}^\ell \Delta_{n,b}^{(1)}(\ell-j+1) PG(X_{j-1})\\
+\sum_{\ell=3}^nPG(X_{\ell-1}) \sum_{j=1}^{\ell-2}\Delta_{n,b}^{(2)}(\ell-j)P G(X_{j-1})\\
+PG(X_0)\left\{\sum_{\ell=1}^nw_{n,b}(\ell)Q_\ell + \sum_{\ell=1}^n \Delta_{n,b}^{(1)}(\ell) PG(X_{\ell-1})\right\}\\
-PG(X_n)\left\{\sum_{j=1}^nw_{n,b}(n-j)Q_j-\sum_{j=1}^n\Delta_{n,b}^{(1)}(\ell-j)PG(X_j)\right\}\\
+\Delta_{n,b}^{(2)}(0)\sum_{\ell=1}^n\left(PG(X_{\ell-1})\right)^2 -\left(w_{n,b}(0)-\Delta_{n,b}^{(2)}(1)\right)\sum_{\ell=1}^nPG(X_\ell)Q_\ell\\
-\Delta_{n,b}^{(2)}(1)\left(PG(X_n)\right)^2 -w_{n,b}(n-1)PG(X_n)PG(X_0).\end{multline*}

Using A\ref{A4}, (\ref{momentV4}), the martingale-difference property of $\{Q_\ell,\;\ell\geq 1\}$, the smoothness of $w_b$,  we derive that for $p>1$ as in (\ref{momentV4}),
\begin{equation}\label{boundzeta2}\PE^{1/p}\left(|\zeta_n|^p\right)\leq c\, c_n^{-1+\frac{1}{2}\vee\frac{1}{p}},\;\;\;n\geq 3,\end{equation}
for some finite constant $c$. 

By martingale approximation for linear partial sums (see e.g. the proof of Proposition \ref{llnMC} below), for any sequence of real numbers $\{a_{n,\ell},\; 1\leq \ell\leq n\}$,
\begin{multline}\label{proofthm3eq1}
\sum_{\ell=1}^na_{n,\ell} h(X_\ell)=\sum_{\ell=1}^n a_{n,\ell}Q_\ell + \epsilon_{n,1},\;\\
\mbox{ where } \;\; \PE\left(\left|\sum_{\ell=1}^n a_{n,\ell}Q_\ell\right|^\alpha\right)\leq c\left(\sum_{\ell=1}^n|a_{n,\ell}|^{\alpha\wedge 2}\right)^{1\vee\frac{\alpha}{2}},\\
\;\mbox{ and } \;\;\PE\left(|\epsilon_{n,1}|^\alpha\right)\leq c\left(|a_{n,1}|+|a_{n,n}|+\sum_{\ell=2}^n|a_{n,\ell}-a_{n,\ell-1}|\right)^\alpha,\end{multline}
provided $\sup_{n\geq 0}\PE\left(V_2^\alpha(X_n)\right)<\infty$. We use (\ref{proofthm3eq1}) to bound the term $R_n$ as given in (\ref{termRn}) and obtain for all $n\geq 1$:
\[\PE\left(\left|R_n\right|^p\right)\leq c n^{-p}c_n^{p},\]
for some finite constant $c$.

By standard martingale inequalities, we obtain the bound
\[\PE\left(\left|\sum_{\ell=1}^n\sum_{j=1}^{\ell-1} w_{n,b}(\ell-j)Q_\ell Q_j\right|^p\right)\leq c \left(\frac{c_n}{n}\right)^{\frac{p}{2}}n^{-\frac{p}{2}+1\vee\frac{p}{2}}.\]

\begin{flushright}
$\square$
\end{flushright}

\subsection{Proof theorem \ref{thm4}}\label{proofthm4}

If $c_n=o(n)$ and $p\geq 2$, then from Theorem \ref{thm3}, $\Gamma_{n,b}^2(h)=n^{-1}\sum_{\ell=1}^nQ_\ell^2+o_P(1)$. Given the ergodicity assumption $\C(1,V_2^2,V_3)$ and $\sup_{k\geq 0}\PE(V_3^q(X_k))<\infty$, it follows from  Proposition \ref{llnMC} that the term $n^{-1}\sum_{\ell=1}^n Q_\ell^2$ converges in probability to the limit
\[\int\mu(dx)\int P(x,dy)\left(G(y)-PG(x)\right)^2,\]
which is easily seen to be equal to $\sigma^2(h)$. This proves the first part of the theorem.

From now on, we assume that $c_n=n$. Define $W_{n,\ell}=\frac{Q_\ell}{\sqrt{n}\sigma(h)}$. Then by (\ref{proofthm3eq1}) with $a_{n,\ell}\equiv \frac{1}{\sqrt{n}\sigma(h)}$, 
\[\sum_{\ell=1}^n\frac{h(X_\ell)}{\sqrt{n}\sigma(h)}=\sum_{\ell=1}^n W_{n,\ell}+\epsilon_{n,2},\;\; \mbox{ where }\;\;\epsilon_{n,2} \mbox{ converges in probability to zero}.\]
Define $\lfloor x\rfloor$ as the largest integer smaller or equal to $x$ and for $0\leq t\leq 1$, we introduce
\[B_n(t)=\sum_{\ell=1}^{\lfloor nt\rfloor}W_{n,\ell},\;\mbox{ and } \;Z_n(t)=\int_0^tw_b(t-u)dB_n(u).\]
Since $B_n$ has jumps only at times $\ell/n=\ell/c_n$, we see that $Z_n(\ell c_n^{-1})=\sum_{j=0}^{\ell-1} w_b((\ell-j)c_n^{-1})W_{n,j+1}$. It is also easy to see that the term $R_n$ in (\ref{termRn}) can be written as 
\begin{multline*}
R_n=2B_n^2(1)\int_0^1 (1-u)w_b(u)du \\
-2B_n(1)\int_0^1\left(\int_0^tw_b(u)du +\int_0^{1-t} w_b(u)du\right)dB_n(t) + \epsilon_{n,3},\end{multline*}
where $\epsilon_{n,3}$ converges in probability to zero. Thus
\begin{multline*}
\Gamma_{n,b}^2(h)=n^{-1}\sum_{\ell=1}^n Q_\ell^2+2\sigma^2(h)\sum_{\ell=1}^n\frac{Q_\ell}{\sigma(h)\sqrt{n}}\sum_{j=1}^{\ell-1} w_b\left(\frac{\ell-j}{c_n}\right)\frac{Q_j}{\sigma(h)\sqrt{n}} + R_n + \zeta_n\\
=\sigma^2(h)\sum_{\ell=1}^nW_{n,\ell}^2+2\sigma^2(h)\sum_{\ell=1}^nW_{n,\ell} Z_n\left((\ell-1)c_n^{-1}\right) + R_n + \zeta_n\;\;\;\\
=\sigma^2(h)\sum_{\ell=1}^nW_{n,\ell}^2+2\sigma^2(h)\int_0^1Z_n(t)dB_n(t) + 2B_n^2(1)\sigma^2(h)\int_0^1 w_b(u)(1-u)du \\
-2B_n(1)\sigma^2(h)\int_0^1g_b(u)dB_n(t) +\epsilon_{n,4},\end{multline*}
where $g_b(u)=\int_0^tw_b(u)du +\int_0^{1-t} w_b(u)du$, and $\epsilon_{n,4}$ converges in probability to zero. 

From the assumptions, $\sup_{\ell\geq 0}\PE\left(|Q_\ell|^{2+\epsilon}\right)<\infty$, for some $\epsilon>0$. Therefore, by the functional central limit theorem for martingales, $B_n\stackrel{w}{\to} B$, where $B=\{B(t),\;0\leq t\leq 1\}$ is the standard Brownian motion. 
By the continuous mapping theorem, $(B_n,Z_n)\stackrel{w}{\to} (B, Z)$, where $Z(t)=\int_0^t w(t-u)dB(u)$. And by the weak convergence of stochastic integrals (see, e.g., Theorem 2.2 in \cite{kurtzetprotter91}), \[\left\{\left(B_n(t),\int_0^tZ_n(t)dB_n(t),\int_0^1 g(u)dB_n(u), \sum_{\ell=1}^n W_{n,\ell}^2\right), 0\leq t\leq 1\right\}\] converges weakly to the stochastic process \[\left\{\left(B(t),\int_0^tZ(u)dB(u),\int_0^1 g(u)dB(u), 1\right),\;0\leq t\leq 1\right\}.\] As the remainders $(\epsilon_{n,2},\epsilon_{n,4})$ converges in probability to $0$, this entails that
$\left(\sum_{\ell=1}^n\frac{h(X_\ell)}{\sqrt{n}\sigma(h)},\Gamma_{n,b}^2(h)\right)$ converges weakly to the limit 
\begin{multline*}
\left(B(1),\right.\\
\left.\sigma^2(h)\left(1+2\int_0^1Z(u)dB(u) + 2B^2(1)\int_0^1 (1-u)w(u)du -2B(1)\int_0^1g(u)dB(t)\right)\right).\end{multline*}
 The conclusion of the theorem follows by the continuous mapping theorem.

\begin{flushright}
$\square$
\end{flushright}

\subsection{Proof Theorem \ref{theoclt1}}\label{prooftheoclt1}
By (\ref{Hoeffdecomp}), we have:
\[\sigma_{n}^{-1}\left[U_n(h_n)-\theta_n{n\choose 2} \right] = \sigma_n^{-1}(n-1)\sum_{\ell=1}^nh_n(X_\ell) +\sigma_{n}^{-1} \sum_{\ell=2}^n\sum_{j=1}^{\ell-1}Q_{n,\ell,j}+\sigma_n^{-1}\zeta_n.\]
(\ref{boundzeta}) gives 
\[\sigma_n^{-2}\PE\left(|\zeta_n|^2\right)\leq cn^{-1}\sigma_{n,1}^{-2}\nnorm{\bar h_{n,2}}^2_{2,\barV_2}=o(1)\]
by (\ref{condeq1}). This  shows that $\sigma_n^{-1}\zeta_n$ converges in probability to zero.

Now, by the  martingale property, we have
\begin{multline*}
\PE\left[\left(\sum_{\ell=2}^n\sum_{j=1}^{\ell-1}Q_{n,\ell,j}\right)^2\right]=\sum_{\ell=2}^n \PE\left[\left(\sum_{j=1}^{\ell-1}Q_{n,\ell,j}\right)^2\right]=\sum_{\ell=2}^n\sum_{j=1}^{\ell-1}\PE\left(Q_{n,\ell,j}^2\right)\\
+2\sum_{\ell=2}^n\sum_{k=1}^{\ell-2}\sum_{j=k+1}^{\ell-1}\PE\left(Q_{n,\ell,j}Q_{n,\ell,k}\right)=o(\sigma_n^2),\end{multline*}
by (\ref{momentQ}), Lemma \ref{lem:techlem1} and (\ref{condeq1}).

\begin{flushright}
$\square$
\end{flushright}

\begin{lemma}\label{lem:techlem1}
Under the assumptions of Theorem \ref{theoclt1}, 
\[\left|\PE\left(\sum_{k=1}^{\ell-2}\sum_{j=k+1}^{\ell-1}Q_{n,\ell,j}Q_{n,\ell,k}\right)\right|\leq cn\nnorm{\bar h_{n,2}}_{2,\barV_2}^2,\;\;3\leq \ell\leq n.\]
\end{lemma}
\begin{proof}
Fix $1\leq k<\ell$ and define 
\[T_k=T_{n,\ell,k}:=\PE\left(\sum_{j=k+1}^{\ell-1}Q_{n,\ell,j}Q_{n,\ell,k}\vert \F_{k-1}\right),\]
so that 
\[\left|\PE\left(\sum_{k=1}^{\ell-2}\sum_{j=k+1}^{\ell-1}Q_{n,\ell,j}Q_{n,\ell,k}\right)\right|\leq \sum_{k=1}^{\ell-2} \PE\left(|T_k|\right).\]
For $m\geq 0$, define
\begin{multline*} \Upsilon_{2,m}(x_{j-1},x_{k-1},x_k)= \int P(x_{j-1},dx_j) \int P^{m}(x_j,dx_{\ell-1}) \\
\times\int P(x_{\ell-1},dx_\ell) \Lambda_{n,2}(x_{j-1},x_{\ell-1};x_j,x_\ell)\Lambda_{n,2}(x_{k-1},x_{\ell-1};x_k,x_\ell),\end{multline*}
\[\Upsilon_{1,m}(x_{k-1},x_k):=\int\left\{P^{m}(x_k,dx_{j-1})-\mu(dx_{j-1})\right\}\Upsilon_{2,\ell-m-k-2}(x_{j-1},x_k,x_{k-1}).\]
Then almost surely we have:
\[T_k=\PE\left(\sum_{j=0}^{\ell-k-1}\Upsilon_{1,j}(X_k,X_{k-1})\vert\F_{k-1}\right).\]
The bound (\ref{boundQ}) and the Cauchy-Schwartz inequality imply that
\begin{multline*}
\left|\int P(x_{\ell-1},dx_\ell)\Lambda_{n,2}(x_{j-1},x_{\ell-1};x_j,x_\ell)\Lambda_{n,2}(x_{k-1},x_{\ell-1};x_k,x_\ell)\right|\\
\leq c\,\nnorm{\bar h_{n,2}}^2_{2,\barV_2}\left(\barW(x_{\ell-1},x_j)+\barW(x_{\ell-1},x_{j-1})\right)\left(\barW(x_{\ell-1},x_k)+\barW(x_{\ell-1},x_{k-1})\right).\end{multline*}
We combine this with B\ref{B2} to conclude that for all $m\geq 0$, 
\[\left| \Upsilon_{2,m}(x_{j-1},x_{k-1},x_k)\right|\leq c\,\nnorm{\bar h_{n,2}}^2_{2,\barV_2}\,\overline{\mathcal{U}}_1(x_{j-1})\overline{\mathcal{U}}_2(x_k,x_{k-1}). \]
By the short-range dependence assumption $\C(1,\overline{\mathcal{U}}_1,\overline{\mathcal{V}}_1)$, it follows that for any $n\geq 0$,
\begin{multline*}
\left|\sum_{j=0}^{\ell-k-1}\Upsilon_{1,j}(X_k,X_{k-1})\right|\\
\leq\sum_{j=0}^{\ell-k-1} \left|\int\left\{P^j(x_k,dx_{j-1})-\mu(dx_{j-1})\right\}\Upsilon_{2,\ell-j-k-2}(x_{j-1},x_{k-1},x_k)\right|\\
\leq c\,\nnorm{\bar h_{n,2}}^2_{2,\barV_2} \overline{\mathcal{V}}_1(x_k)\overline{\mathcal{U}}_2(x_k,x_{k-1}),\end{multline*}
for some finite constant $c$. We conclude that
\[|T_k|\leq c\,\nnorm{\bar h_{n,2}}^2_{2,\barV_2} \PE\left(\overline{\mathcal{V}}_1(X_k)\overline{\mathcal{U}}_2(X_k,X_{k-1})\vert \F_{k-1}\right).\]
The lemma follows.
\end{proof}

\section{Appendix A: A weak law of large numbers for Markov chains}
\begin{prop}\label{llnMC}
Let $\{X_n,\;n\geq 0\}$ be a Markov chain with invariant distribution $\mu$ and transition kernel $P$. Suppose that there exist measurable functions $V_1\leq V_2:\;\X\to [1,\infty)$ such that 
\begin{equation}\label{rate2}\sum_{k\geq 0}\|P^k(x,\cdot)-\mu\|_{V_1}\leq cV_2(x),\;x\in\X,\end{equation}
for some finite constant $c$. Suppose also that $\v_n:=\PE(V_2^p(X_n))<\infty$ for each $n\geq 0$ and for some $p\in (1,2]$. Let $\{f_n,\;n\geq 1\}$ be such that $f_n, Pf_n\in\L_{V_1}$ and let $\{a_{n,k}, 0\leq k\leq n\}$ be a sequence of real numbers such that 
\begin{multline*}
\nnorm{f_n}_{p,V_1}^p\left(\sum_{k=1}^n|a_{n,k}|\right)^{-p}\sum_{k=1}^n|a_{n,k}|^p \v_k^p\to 0,\;\;\\
\;\;\mbox{ and }\;\;|Pf_n|_{V_1}\left(\sum_{k=1}^n|a_{n,k}|\right)^{-1}\sum_{k=1}^n|a_{n,k}-a_{n,k-1}|\v_{k-1}\to 0.\end{multline*}
 Then, as $n\to\infty$, $\left(\sum_{k=1}^n|a_{n,k}|\right)^{-1}\sum_{k=1}^na_{n,k}\left(f_n(X_k)-\mu(f_n)\right)$ converges in probability to zero.
\end{prop}
\begin{proof}
Define $S_n=\sum_{k=1}^na_{n,k}\left(f_n(X_k)-\mu(f_n)\right)$ and $g_n(x)=\sum_{j\geq 0}(P^jf_n(x)-\mu(f_n))$. Under (\ref{rate2}), $|g_n(x)|\leq c|f_n|_{V_1}V_2(x)$ and $|Pg_n(x)|\leq c|Pf_n|_{V_1}V_2(x)$. By the Poisson equation, $f_n(x)-\mu(f_n)=g_n(x)-P g_n(x)$ which implies that
\begin{multline*}
S_n=\sum_{k=1}^na_{n,k}\left(g_n(X_k)-P g_n(X_{k-1})\right)+\sum_{k=1}^n(a_{n,k}-a_{n,k-1})Pg_n(X_{k-1})\\\
+\left(a_{n,0}Pg_n(X_0)-a_{n,n}Pg_n(X_n)\right).\end{multline*}
where the martingale array $\sum_{k=1}^na_{n,k}\left(g_n(X_k)-P g_n(X_{k-1})\right)$ satisfies
\[\PE\left(|\sum_{k=1}^na_{n,k}\left(g_n(X_k)-P g_n(X_{k-1})\right)|^p\right)\leq C\nnorm{f_n}_{p,V_1}^p\sum_{k=1}^n|a_{n,k}|^p\v_k^p.\]
The last inequality follows by noting that $g_n(x)-P g_n(y)=f_n(x)-\mu(f_n)-Pg_n(x)-P g_n(y)$ and by conditioning on $\F_{k-1}$. Thus, under the stated assumptions, $\left(\sum_{k=1}^n|a_{n,k}|\right)^{-1}S_n$ converges in probability to zero.
\end{proof}
\begin{rem}\label{remlln}
An important special case is the case where $a_{n,\ell}=1$ and $\sup_{n\geq 0}\PE\left(V_2^p(X_n)\right)<\infty$. In this case it is enough to have $n^{-1+1/p}\nnorm{f_n}_{p,V_1}\to 0$. If in addition it is true that $\sup_{x\in\X}PV_1^p(x)/V_1^p(x)<\infty$, then clearly $\nnorm{f_n}_{p,V_1}\leq c |f_n|_{V_1}$ and the law of large number holds if $n^{-1+1/p}|f_n|_{V_1}\to 0$.
\end{rem}

\end{document}